\documentclass[12pt,aps,floatfix]{revtex4}

\usepackage{graphicx}
\usepackage{amssymb}
\usepackage{amsmath}
\usepackage{amsfonts}
\usepackage{epstopdf}
\usepackage{epsfig}
\usepackage{wrapfig}

\newtheorem{theorem}{Theorem}

\newtheorem{definition}{Definition}

\begin{document}

\title{{\LARGE An Introduction To Geometric Prequantization}}
\author{Joseph Geraci}
\affiliation{MaRS - TMDT \\ Signaling Biology}
\affiliation{University of Toronto \\ Mathematics}
\maketitle

\section*{\protect\underline{{\protect\Large Abstract}}}


Classical mechanics has a natural mathematical setting in symplectic geometry and it may be asked if the same is true for quantum mechanics. More precisely, is it possible to capture certain quantum idiosyncrasies within the symplectic framework of classical mechanics?  Beginning with Dirac, the idea of quantizing a classical system involved associating the phase space variables with Hermitian operators which act on some Hilbert space, as well as associating the Poisson bracket on the phase space with the commutator for the operators. Mathematically the phase space is associated with some symplectic manifold $M$ and the non-degenerate 2-form $\omega$, which comes with $M$. Geometric prequantization is a process by which one does this in a mathematically ``rigorous'' manner and we shall attempt to just introduce the methods here.  We do this by exploring this contruction for ($\mathbb{R}^{2n}, \sum_{i=1}^n dp_i \wedge dq_i$) which avoids using complex line bundles in any non-trivial way. One should note however that the Hilbert space one obtains is in fact too ``big'', in the sense that it has too many functions in order to correspond with actual physically significant Hilbert spaces. Geometric quantization remedies this situation but it should be noted that not all manifolds are prequantizable. We shall not discuss either of these issues however. 

\newpage

\section{Some Classical Mechanics}

Classical mechanics has a very rich and beautiful mathematical structure and symplectic geometry is the natural mathematical setting to study such systems. Let us recall the Hamiltonian description of a classical system. The simple example of a particle in three dimensions, $\mathbb{R}^3$, trapped in a potential $V(q_1,q_2,q_3)$ will suffice. Note that this three dimensional system has a six dimensional phase space with variables $(q_i,p_i)$ where $i=1,2,3$ and the $q_i$ are the position variables and the $p_i$ are the momenta. Let the Hamiltonian function H, which represents the energy of the system, be given by 

\[ H = \frac{1}{2m}(p_1^2 + p_2^2 + p_3^2) + V(q_1,q_2,q_3). \] 

Now, the phase variables are functions of time and the evolution of the state of the system is given by the solutions to Hamilton's equations

\begin{eqnarray}  \frac{dq_i}{dt} &=& \frac{\partial H}{\partial p_i} \label{Ham1} \\
\frac{dp_i}{dt} &=& -\frac{\partial H}{\partial q_i}.\label{Ham2} \end{eqnarray}

The solutions to these equations represent the sate of the particle in the phase space and not the configuration space. What this means is that one obtains information about the position and momentum of the particle for any given time. The jump to analyze a particle in any dimension $\mathbb{R}^n$, is immediate and one should note that the phase space will now be $\mathbb{R}^{2n}$. Finally, we must mention a very important algebraic structure that comes with the phase space known as the Poisson bracket. It is a bilinear, anti-symmetric operation on pairs of $C^\infty$ functions (which act on the phase space) given by 

\[ \{ f,g \} = \sum_{i=1}^n \left(\frac{\partial f}{\partial q_i} \frac{\partial g}{\partial p_i} - \frac{\partial f}{\partial p_i} \frac{\partial g}{\partial q_i} \right ). \]

The geometry comes in when one asks if we could have a coordinate free formulation of Hamiltonian Mechanics so that we could consider more interesting phase spaces. Notice that our phase space above was even dimensional. Consequently, if we restrict ourselves to even dimensional manifolds there is a beautiful way to extend this analysis. Specifically, we will use symplectic manifolds to accomplish this. Recall that a symplectic manifold M is an even dimensional manifold equipped with a non-degenerate and closed 2-form $\omega$. Our friend $\mathbb{R}^{2n}$ equipped with the 2-form $\omega_o = \sum_{i=1}^n dp_i \wedge dq_i$ is a symplectic manifold for example. A rich source of symplectic manifolds are given by the cotangent bundles of manifolds since these are always even dimensional and they come equipped with a canonical 2-form, in general. 

Now, recall that a symplectomorphism is a diffeomorphism $T\!:M \rightarrow M$ which preserves the symplectic structure i.e. $T^* \omega = \omega$. Thanks to the Darboux Theorem we have that any symplectic manifold is locally symplectomorphic to $\left ( \mathbb{R}^{2n}, \omega_o \right )$. This will be pertinent in a little while. For now let's get back to Hamilton's equations. Here we go from our Hamiltonian function $H$ to the Hamiltonian vector field  $X_H$ by using the symplectic 2-form $\omega$  by considering the map 

\begin{equation} \tilde{\omega}: TM \rightarrow T^* M \label{map} \end{equation} given by \[ X_H = \tilde{\omega}^{-1}(dH).  \]

Equivalently \begin{equation} \omega(X_H, \cdot) = i_{X_H} \omega = -dH. \label{symp} \end{equation}

(Note that the negative appearing in the above equation is just a matter of common convention.)

So now we must see what all this means. The first thing to notice is that since $\omega$ is non-degenerate, equation (\ref{symp}) has a unique solution for the vector fields $X_H$.   Let's go back to $(\mathbb{R}^{2n}, \omega_o)$ and work out some details. First we need dH which is given by 

\[ dH = \sum_{i=1}^n \left ( \frac{\partial H}{\partial q_i} dq_i + \frac{\partial H}{\partial p_i} dp_i. \right ) \]

Then the equation \begin{equation} i_{X_H} \omega_o  = -dH \label{sympR} \end{equation} implies that \begin{equation} X_H = \sum_{i=1}^n \left ( \frac{\partial H}{\partial p_i} \frac{\partial}{\partial q_i} - \frac{\partial H}{\partial q_i} \frac{\partial}{\partial p_i}\right ). \label{hvec} \end{equation}

This means that the solutions to Hamilton's equations given by equations ($\ref{Ham1}$) and ($\ref{Ham2}$) are the flow of $X_H$ given above. Therefore, because of the Darboux Theorem, we know that this holds locally on any symplectic manifold. In fact, equation ($\ref{symp}$) holds globally on any symplectic manifold and the map ($\ref{map}$) gives an isomorphism between the tangent and cotangent spaces at each point $m \in M$. So we have a natural way of identifying the space of complex vector fields on M with the space of complex 1-forms. In this way, the vector fields which correspond to exact 1-forms are referred to as global Hamiltonian vector fields and those corresponding to closed 1-form are referred to as local Hamiltonian vector fields \cite{wood}. A very important property of $X_H$ (and one which can be taken as a definition of a Hamiltonian vector field) is that the Lie derivative of $\omega$ with respect to $X_H$ is zero or $\mathcal{L}_{X_H} \omega = 0$.  This is obvious if you use Cartan's magic formula. One finds

\[ \mathcal{L}_{X_H} \omega = i_{X_H}(d\omega) + d(i_{X_H}(\omega) = -d(dH) = 0 \] since $\omega$ is closed by definition. One last important thing to notice is that equation ($\ref{hvec}$) immediately implies the identity

\[ \{ \phi, \psi \} = \omega(X_\phi,X_\psi) = X_\psi \phi \] which is true for any symplectic manifold, where $\phi$ and $\psi$ are in $C^{\infty}(M)$.
 
 So physically, what is going on? One should imagine that the the points on the manifold M are in 1-1 correspondence with the states of the physical system one is considering.  Further, say that the system is in a state $m \in M$ at some time $t_o$. Then the system will be in the state $m' = \phi_{t'-t_o}(m)$ where $\phi_t$ is the flow of the Hamiltonian vector field $X_H$.  Further, the Lie algebra $(C^{\infty}(M), \{\cdot, \cdot\})$ with plays a fundamental role in what we wish to do in this work. The observables in classical mechanics (or the measurable quantities) are represented by functions in $C^{\infty}(M)$. Now, Hamilton's equations may be rewritten using the Poisson bracket i.e. we may write 
 
 \[ \frac{d}{dt}q_i = \{H,q_i \}, \hspace{2cm}  \frac{d}{dt}p_i = \{H,p_i \} \] and in general one finds that the rate of change of an observable $f \in C^{\infty}(M)$ is given by the Poisson bracket of $f$ with $H$ i.e. 
 \[ \frac{d}{dt} f = \{ H, f \}. \] A very important observation then is to note that $f$ is a conserved quantity (that is, the value of $f$ evaluated at some state of the system at some time $t$ does not depend on $t$) if $\{H,f\} = 0$.

Now, if one agrees that $X_H$ should depend linearly on the rate of change of the Hamiltonian, i.e. the energy of the system, then we need a way to associate $X_H$ with $dH$. But $M$ comes equipped with $\omega$ where we can write $\omega(X_H, V) = -dH(V)$ for any vector field $V$. The dynamics of our physical system are just the flow lines of $X_H$ as mentioned above and the conservation of energy, which is an essential ingredient of a Hamiltonian system, is then given by $dH(X_H) = 0$. This last equation just expresses the physical fact that the energy is a constant for Hamiltonian systems along the flow lines of $X_H$.

\section{Some Quantum Mechanics}

Mathematically, what is essential to Quantum Mechanics(QM) is that we have a complex Hilbert space $\mathcal{H}$, some Hamiltonian (self-adjoint operator) $\mathbb{H}$ on $\mathcal{H}$, and of course the Schrodinger equation \[ \mathbb{H} \psi = -i \hbar \partial_t \psi .\] The solutions to this equation are elements in $\mathcal{H}$ and represent states of the quantum system with Hamiltonian $\mathbb{H}$. The observables here are now self-adjoint operators and the result of a measurement is an eigenvalue of some self-adjoint operator. The outcome of a measurement, i.e. which eigenvalue occurs, is in general not completely predictable and therefore there is a probabilistic interpretation to QM.  What is essential for us to understand is that the state $\psi$ and $c\psi$ for $c\neq 0 \in \mathbb{C}$ are physically identical. This means that all non-zero points of the complex line $\{c\psi | c \in \mathbb{C} \}$ represent the same state. This is included in the postulates of QM and basically it means that the "direction" and not the magnitude of the state is what matters. What this means geometrically is that the states of a QM system form a projective Hilbert space $\mathbb{P} \mathcal{H}$. For example, the finite dimensional Hilbert space most commonly found in quantum computation is the complex plane $\mathbb{C}$ and so one often speaks of the states of the quantum computer being in $S^2$ which is homeomorphic to $\mathbb{CP}^1$.  So to continue, we now understand that a point in $\mathcal{H}$ does not determine a state of a QM system uniquely, however a point in  $\mathbb{P} \mathcal{H}$ does. The problem here is that $\mathbb{P} \mathcal{H}$ is quite difficult to work with, not to mention it is infinite dimensional in general. Instead we shall briefly discuss the situation that occurs if we just study the unit sphere in $\mathcal{H}$, $S \mathcal{H} = \{ \psi \in \mathcal{H} | \langle \psi, \psi \rangle = 1 \}$. Though this space is also in general infinite dimensional, it will give us the ability to analyze the dynamics of a QM system in such a way that it will enable us to produce an analogy with classical mechanics where certain phenomenon which are thought to only occur in a QM system, also appear in the classical framework.

By studying $S \mathcal{H}$ instead of the projective space, we gain in simplicity but we loose uniqueness. Two points in  $S \mathcal{H}$ which differ by a multiple of unit modulus represent the same state, i.e. $\psi$ and $e^{i\theta} \psi $ both correspond to the same QM state. Now, let $\pi$ be the projection 

\begin{equation} 
\pi : S \mathcal{H} \rightarrow \mathbb{P}\mathcal{H}. \label{proj} \end{equation} 

The inverse image of a point in $\mathbb{P} \mathcal{H}$ is now a circle in $S \mathcal{H}$ and two points on this circle differ by a phase factor $e^{i\theta}$ which often shows up in the beginning of any course on QM. 

Finally, we need to see what the quantum dynamics are when our model consists of $S \mathcal{H}$. First, if $\mathbb{H}$ is the Hamiltonian of our system then we can define the vector field $X_H$ on $\mathcal{H}$ via 

\begin{equation} X_H = \frac{i}{\hbar} \mathbb{H} \psi. \end{equation} This allows us to realize the Schrodinger equation as the flow of $X_H$ i.e. $\frac{d\psi}{dt} = X_H(\psi).$ We have the following theorem.

\begin{theorem} Given that $\mathbb{H}$ is self-adjoint, $X_H$ is tangent to the unit sphere $S \mathcal{H}$.
\end{theorem}

The way to see this is to realize that the equation of the flow of $X_H$ given above implies 
\[ \dot{\psi}  = \frac{i}{\hbar} \mathbb{H} \psi. \] But $\mathbb{H}$ is self-adjoint which means that the propagator for this differential equation is a unitary operator. More specifically we have the propagator $U(t)$ where

\[ U(t) = exp(\frac{it}{\hbar} \mathbb{H}). \] This is the flow for $X_H$ and it is obtained by just solving formally for $\psi$. But note that $U(t)$, being unitary, preserves length and it can therefore be shown to be a one parameter group of diffeomorphism of $S \mathcal{H}$. By definition then, the vector field will be tangent to $S \mathcal{H}$. 

Now that we have the above observations, we know that quantum dynamics can be modeled mathematically by using $S \mathcal{H}$ and $X_H$ where we now have uniqueness up to a phase factor. Can we obtain some analogous situation by using the classical mechanical framework presented earlier? The answer to this question is an affirmative and it relies on prequantization. 

\section{What is Prequantization?}

Let us give a formal definition here so that we can have some mathematical footing. 

\begin{definition} Let $\mathcal{P}$ be a sub-algebra of the Poisson algebra $(C^{\infty}(M), \{\cdot, \cdot \})$
which contains the constant function 1. A prequantization of $\mathcal{P}$ is a linear map $\Omega$ from $\mathcal{P}$ to the linear space of symmetric operators which leave fixed some dense domain within some separable Hilbert space. The following conditions must be met:
\begin{enumerate}

\item $\Omega( \{f,g\}) = \frac{i}{\hbar}[\Omega(f), \Omega(g)]$ where $f,g \in C^{\infty}(M)$
\item $\Omega(1) = I$ 
\item When $X_H$ is complete then $\Omega$ is self-adjoint. 
\end{enumerate}

\end{definition}

So what this means is that prequantization is a procedure to construct an isomorphism from the Poisson algebra of our given manifold to the space of linear self-adjoint operators which act on some Hilbert space in such a way so that the commutation relation $(1)$ in the definition is satisfied. Where does the Hilbert space come from and what is the nature of the map $\Omega$? We shall attempt to answer these questions briefly and then turn to the prequantization of $(\mathbb{R}^{2n}, \omega_o)$.  

So we need a Hilbert space and all we have is the phase space $(M,\omega)$. There are several ways of obtaining a Hilbert space in this situation in general and we take the most naive way here. (Refer to \cite{puta} for more on this.) This involves introducing a complex line bundle $L$ over $M$. One then takes all the sections of $L$, $\Gamma(L)$, as the Hilbert space. However, the Hilbert space obtained in this way is too large and is cut down by geometric quantization, but it is a step in the right direction. Thus, what we do is ``attach'' an operator $\hat{f}: \Gamma(L) \rightarrow \Gamma(L)$ to each classical observer $f \in C^{\infty}(M)$. Specifically, $\hat{f}(s) = \nabla_{X_f} s + 2\pi i f s$ where, as before, $X_f$ is the Hamiltonian vector field generated by $f$ and $\nabla$ is a covariant derivative on $L$. $L$ is then called the prequantization bundle of $M$ if $\{\hat{f},\hat{g}\} = \hat{f}\hat{g} - \hat{g} \hat{f} = [\hat{f},\hat{g}]$ as required by the above definition. Just as a last note, this requirement can be shown to be identical to the existence of a covariant derivative $\nabla$ whose curvature is $\omega$. This only occurs for manifolds with integral cohomology, and therefore it is clear that not all manifolds can be prequantized.  

The map $\Omega$ which makes the association $f \rightarrow \hat{f}$ actually comes from the Lie algebra morphism $f \rightarrow \frac{h}{i} X_f$. However, requirement (2) from the definition above is not met since this map is not 1-1 being that the kernel of this map consists of all constant functions on $M$. We shall see in the next section that this will be remedied by the introduction of a complex line bundle. For now, we should be convinced that we indeed have this morphism of Lie algebras. Let $V(M)$ be the set of globally defined Hamiltonian vector fields. The commutator bracket gives $V(M)$ the structure of a Lie algebra. 

\begin{theorem} The map $f \rightarrow X_f$ is a homomorphism of Lie algebras from $(C^{\infty}(M),\{ \cdot , \cdot \})$ to $(V(M), [ \cdot, \cdot ])$. 
\end{theorem}

$\bf{Proof:}$ We have to show that this map is linear over $\mathbb{C}$ and that it preserves brackets. First of all, we have 
\[ \{f,g\} = X_f(g) = -i_{X_f}(i_{X_g} \omega) = \omega(X_f,X_g), \]
and therefore we have bilinearity over $\mathbb{C}$ as well as skew symmetry. Now we argue as follows. If $X_f$ and $X_g$ are the Hamiltonian vector fields of $f$ and $g$ then $[X_f, X_g]$ is the Hamiltonian vector field of $\omega(X_f,X_g)$ or more precisely we have

\[ i_{[X_f,X_g]}\omega = -d\omega(X_f,X_g) \] which follows directly from equation ($\ref{symp}$).  Continuing in this way we have 

\begin{eqnarray*} -d\omega(X_f, X_g) &=& i_{X_{\omega(X_f,X_g)}} \omega \\
&=& i_{X_{\{f,g\}}} \omega \end{eqnarray*}

which means that we have \[ [X_f, X_g] = X_{\{f,g\}}. \] This means that this morphism indeed does preserve brackets if we take as given that the Poisson bracket satisfies the Jacobi identity. $\blacksquare$

\section{Prequantization of $(\mathbb{R}^{2n}, \omega_o)$}

Here we finally arrive at a simple demonstration of prequantization. Our symplectic manifold is $M=\mathbb{R}^{2n}$ which is the phase space of a particle in n dimensional euclidean space. We will demonstrate that the dynamical situation described above for quantum mechanics, namely the uniqueness up to a phase factor, can be duplicated in classical mechanics. This is accomplished in general by geometric prequantization. Recall specifically we had the projection $\pi :S\mathcal{H} \rightarrow \mathbb{P} \mathcal{H}$ where the inverse image of a point is a circle of phase factors. Note that the circle of phase factors can be represented by the group $U(1)$. Prequantization will also give us another surprise when we study the integral curves of $X_H$. 

Now recall that the canonical 2-form here is $\omega_o = \sum_k dp_k \wedge dr_k$ where our phase space variable are $(\bf{r},\bf{p})$. By introducing a complex line bundle $L$ we will accomplish two things. We will trivially duplicate the freedom in the phase factor found in QM, and also we will be able to extend the Lie algebra homomorphism described above to an injective map. This is a very important step, for we will have the beginning of the goal of canonical quantization originally thought of by Dirac i.e. to have a 1-1 correspondence between classical observables and symmetric operators. We take as our line bundle

\begin{equation} L = \mathbb{R}^{2n} \times U(1) \hspace{2cm} \mathrm{with} \hspace{1cm} \pi: L \rightarrow \mathbb{R}^{2n}. \end{equation}
Specifically $\pi: (\bf{r},\bf{p}, e^{i\theta}) \rightarrow (\bf{r}, \bf{p})$. We can use $\theta$ to parametrise $U(1)$ instead so that our coordinates on $L$ are $(\bf{r},\bf{p}, \theta).$  Now, vector fields on $M$ will be mappings $X:M \rightarrow \mathbb{R}^{2n}$ and on $L$ they will $V:L \rightarrow \mathbb{R}^{2n+1}$. Then we can see from equations (\ref{sympR}) and (\ref{hvec}), that for some $f \in C^{\infty}(M)$
the Hamiltonian Hamiltonian vector field $X_f$ is given by

\[ X_f(\bf{r}, \bf{p}) = \left( \frac{\partial f}{\partial \bf{p}}(\bf{r}, \bf{p}), -\frac{\partial f}{\partial \bf{r}}(\bf{r}, \bf{p}) \right). \]

How about $V_f$? We know that the map $f \rightarrow X_f$ is a Lie algebra homomorphism but that it is not injective. However, we can find a $V_f$ to remedy this situation. We need to select the last coordinate to be a $\mathbb{R}$ -valued function so that each constant function will be distinguished. Let us take 

\[ V_f(\bf{r}, \bf{p},\theta) = \left( \frac{\partial f}{\partial \bf{p}}(\bf{r}, \bf{p}), -\frac{\partial f}{\partial \bf{r}}(\bf{r}, \bf{p}), f(\bf{r},\bf{p}) - \bf{p} \cdot \frac{\partial f}{\partial \bf{p}}(\bf{r},\bf{p}) \right). \]

The mysterious looking term in the $\theta$ coordinate is none other than the negative of the Lagrangian $\Lambda$. Briefly, for the Hamiltonian $H$ we have $\Lambda(\mathbf{r}, \bf{v}) = \mathbf{p} \cdot \frac{\partial H}{\partial \mathbf{p}} - H$. 

This definitely works because now different constant functions will be mapped to different vector fields. So we have an injective Lie algebra homomorphism $f \rightarrow V_f$ as desired as well as uniqueness up to a phase factor that one encounters in QM. 

Previously, it was mentioned that we will get something of a surprise when we studied the integral curves of $X_H$. The truth is that this happens when one looks at $V_H$ instead. For the first two coordinate one obtains the integral curves $\bf{r}(t)$ and $\bf{p}(t)$. For the third we have to solve $\dot{\theta} = - \Lambda$. So

\[  \theta (t) = - \int_{t_o}^t \Lambda ds + \theta (t_o) \] which means we have 
\[ e^{i \theta(t)} = e^{(i\theta(t_o)}  exp\left(-i \int_{t_o}^{t} \bf{p} \cdot \frac{\partial H}{\partial \bf{p}}(\bf{r},\bf{p}) - H(\bf{r}(s),\bf{p}(s)) ds \right) . \]

This is the phase factor in the Feynman path integral. In this approach to QM, a particle is imagined as evolving by taking all paths say from a starting point $a$ and ending at some point $b$. One can ask what the probability is of the particle getting to $b$. Some paths are favored over others and their relative weights are given by the above integral. In classical mechanics there is nothing interesting about this situation. Either a particle gets to $b$ or not, and there is no probabilistic ambiguity. But in the Feynman approach to QM, the probability amplitude is determined by summing over all the paths. We have duplicated a piece of this very quantum mechanical piece of machinery from our simple model. 

What we failed to mention is that $L$ comes with a 1-form $\alpha$ such that $d\alpha=\omega$. This is very important for general considerations. $\alpha$ is given by

\[ \alpha = \bf{p} \cdot d\bf{r} + d\theta \] and so $L$ is a principal $U(1)$ bundle over M with connection $\alpha$ and curvature $\omega$. In this way the vector fields $V_f$ are the unique vector fields  on $L$ (see reference \cite{tuny}) which satisfy 

\[ \pi _* V_f = X_f  \hspace{2cm} \mathrm{and} \hspace{2cm} \alpha(V_f) = f , \] so that the above choice for $V_f$ is not just some ad hoc construction.

A principal fibre bundle $L$ (a fibre bundle with fibre identical to the structure group) with connection $\alpha$, curvature $\omega$ and structure group the circle can be found for most symplectic manifolds. Once one has this much, then an injective Lie algebra morphism from the Poisson algebra to the vector fields on $L$ can be found. The 1-form $\alpha$ comes about during the search for the 1-1 morphism. Again see \cite{tuny} for a fine introduction. \cite{wood} then goes on in much more detail. We should mention that in general the fibre bundle $(L,\alpha)$ is not unique, as it was for us. In general, if M is simply connected then one has uniqueness. For motivation we will mention that a local trivialization can be chosen for $(L, \alpha)$ so that $\alpha$ looks locally like $\alpha = \bf{p} \cdot d\bf{r} + d\theta$ and this corresponds to a local gauge transformation. This is extremely important if one wants this construction to be physically realistic. In fact, the gauge symmetries one finds in QM are considered the basis for electromagnetism. Further, as mentioned previously, the Hilbert space comes about by considering section of $L$. In this case, the sections really are complex valued functions. However, even here, the Hilbert space is too large to be physically useful. Geometric quantization deals with these issues. Please refer to \cite{tuny} and \cite{wood} for details. 

So what did we accomplish? We have the fact that a classical state is a point in $M$ and a circle in $L$ which is analogous to the situation discussed in QM with points in $\mathbb{P}\mathcal{H}$ and circles in $S\mathcal{H}$.  We also obtained a very interesting integral curve that corresponds to the Feynman path integral approach to QM. Most importantly, we have seen that the introduction of the fibre bundle $L$ allowed us to have an injective Lie algebra morphism, which is the basis for associating classical observables with quantum mechanical observables. This is a great beginning to a more complete process known as geometric quantization which is a vast and beautiful subject.

\section{Literature Citations}


\begin{thebibliography}{99}
\bibitem{wood} {D.J. Simms and N.M.J. Woodhouse}, Lecture Notes in Physics - Lectures on Geometric Quantization, Springer-Verlag, (1976)

\bibitem{tuny} {G.M. Tunyman}, What is prequantization, and what is geometric quantization, Proc. Seminar 1989-1990 on Mathematical Structures in Field Theory. CWI Syllabi 39, Math. Centrum, Centrum Wisk. Inform., Amsterdam (1996).

\bibitem{snia} {Jedrzej Sniatycki}, Geometric Quantization and Quantum Mechanics, Applied Mathematical Sciences 30, Springer Verlag,(1980)
(1997).

\bibitem{puta} {Mircea Puta}, Hamiltonian Mechanical Systems and Geometric Quantization, Kluwer Academic Publishers Group (1993).

\bibitem{gotay} {M.J. Gotay,H.B. Grundling and G.M. Tuynman},Obstruction Results in Quantization Theor,J. Nonlinear Sci. 6, 469-498 (1996).
 

\end{thebibliography}
\end{document}